\documentclass[preprint,12pt]{elsarticle}
\usepackage{graphicx}
\usepackage[all]{xy}
\usepackage{amssymb}
\usepackage{amsmath}
\usepackage{amscd}
\usepackage{amsfonts}
\usepackage{pxfonts}

\newcommand{\Q}{\mathbb{Q}}
\newcommand{\R}{\mathbb{R}}
\newcommand{\SNT}{\mathrm{SNT}}
\newcommand{\Z}{\mathbb{Z}}

\newcommand{\p}{_{(p)}}
\newcommand{\Ker}{\mathrm{Ker}}
\newcommand{\Ph}{\mathrm{Ph}}
\numberwithin{equation}{section}

\journal{Topology and its Applications}

\begin{document}

\begin{frontmatter}



\title{On the Gray index conjecture for phantom maps}


\author{Kouyemon Iriye}
\ead{kiriye@mi.s.osakafu-u.ac.jp}
\address{Department of Mathematics and Information Sciences \\
	 Osaka Prefecture University  \\
	 Sakai, Osaka, 599-8531, JAPAN}

\begin{abstract}
We study the Gray index, a numerical invariant for phantom maps. 
It has been conjectured that the only phantom map between finite-type 
spaces with infinite Gray index is the constant map. We disprove this conjecture 
by constructing a counter example.  We also prove that this conjecture is valid 
if the target spaces of the phantom maps are restricted to 
being simply connected finite complexes. 

As a result of the counter example, we can show that $\SNT^{\infty}(X)$ can be 
non-trivial for some space $X$ of finite type.  

\end{abstract}


\end{frontmatter}


\section{Introduction}
\label{}
Throughout this paper all spaces have base-points, and maps and homotopies 
between them are pointed. 

Recall that a map $f:X\to Y$ is called a {\it phantom} map if for any finite {\it dimensional}  
CW-complex $K$ and any map $g:K\to X$, the composition $f\circ g:K\to Y$ is null homotopic.  
Here by a finite dimensional CW-complex we mean one  
with only $n$-cells for n less than some fixed finite number.  

In the literature there is another slightly different notion of a phantom map, 
which may be more common at least when considering stable phantom maps. 
In it a map $h:Z\to W$ is said to be a phantom map if for any finite complex $L$ 
and any map $k:L\to Z$, the composition $h\circ k:L\to W$ is null homotopic, 
e.g., see \cite{Z}. McGibbon gives a full detail of the reason why 
we choose the first definition about a phantom map at the first section of his survey paper 
\cite{Mc}. Needless to say, if $X$ is of finite type then the two definition coincide. 

We write $\Ph(X,Y)$ for the set of homotopy classes of phantom maps from $X$ to $Y$. 
It is, in general, only a pointed set. Its base point is the homotopy class of the constant map. 
If $X$ has the rational homotopy type of a suspension or $Y$ has 
the rational homotopy type of a loop space, 
then $\Ph(X,Y)$ has a natural abelian group structure. 
As it is a huge set unless trivial (e.g., see \cite{Mc} or \cite{R}), 
it is natural to seek for filtrations or invariants 
to distinguish one phantom map from another. 
The Gray index for phantom maps is one of such invariants and goes back to \cite{G}; it 
has been recently studied in \cite{Gh}, \cite{Le}, \cite{LS}, \cite{LS03}, \cite{L} and \cite{MS}. 

If $X$ is a CW-complex and $f:X\to Y$ is a phantom map, then for each natural number $n$ 
there is a map $f_n:X/X_n\to Y$ such that $f\simeq f_n\circ \pi_n$, 
where $X_n$ is the $n$-skeleton of $X$ and $\pi_n:X\to X/X_n$ is the canonical collapsing map. 
The Gray index $G(f)$ for $f$ is the  
least integer $n$ for which $f_n:X/X_n\to Y$ cannot be chosen to be a phantom map. 
If for every $n$ we can choose $f_n$ as a phantom map, then we say that $f$ has infinite 
Gray index, which is denoted by $G(f)=\infty$. 
We denote the set of all homotopy classes of phantom maps 
$f:X\to Y$ with $G(f)\geq n$ by $\Ph^n(X,Y)$. 
The constant map is a phantom map with 
infinite Gray index. It is easy to see that the Gray index is a homotopy 
invariant for phantom maps. Moreover, it does not depend on the choice 
of the CW-structure on $X$ (see \cite{G} or \cite{LS}). 

It is natural to conjecture that every essential phantom map has finite Gray index. Unfortunately, 
McGibbon and Strom \cite{MS} have constructed an essential phantom map 
out of $\mathbb{C}P^{\infty}$ 
with infinite Gray index. The target space in their example is, however, not of finite type. 
This observation has led to the following conjecture. 
\\

\noindent
{\bf Conjecture 1.1.} {\it $\Ph^{\infty}(X,Y)=*$ for finite type domains $X$ and 
finite type targets $Y$.}
\\

Here a space $X$ is called a finite type domain if each of its integral 
homology groups is finitely generated; a space $Y$ is referred to as a finite type 
target if each of its homotopy groups is finitely generated.  
This conjecture is known as the Gray index conjecture. 
But we are able to disprove this conjecture as follows.
\\

\noindent
{\bf Theorem 1.2.} {\it There is a 3-connected space $Y$ of finite type such that \linebreak
$\Ph^{\infty}(\mathbb{C}P^{\infty},\Omega Y)\neq*$.}
\\

For a space $X$, by $\SNT(X)$ we denote the set of homotopy 
types of spaces $Y$ having the same $n$-type for $X$, 
for all $n$. That is, $X^{(n)}$ and $Y^{(n)}$, 
the Postnikov approximations of $X$ and $Y$ for dimension $n$, 
have the same homotopy type for all $n$. 
Ghienne \cite{Gh} introduces a natural filtration on $\SNT(X)$: 
\[
\SNT(X)\supset \SNT^1(X)\supset \cdots \supset \SNT^k(X)\supset\cdots
\supset \SNT^{\infty}(X)=\cap_k\SNT^k(X), 
\]
which has the same algebraic characterization as the Gray index for phantom maps 
and whose precise definition will be given in section 2. 
\\

\noindent
{\bf Corollary 1.3.} {\it $\SNT^{\infty}(\mathbb{C}P^{\infty}\times\Omega^2 Y)\neq*$, 
where $Y$ is the space in Theorem 1.2.}
\\

Although the Gray index conjecture is not true in general, it is valid under additional hypotheses. 
\\

\noindent
{\bf Theorem 1.4.} {\it Let $X$ and $Y$ be nilpotent spaces of finite type. 
If any of the following conditions hold, then $\Ph^{\infty}(X,Y)=*$. 
\begin{enumerate}
\item[(i)]  $X$ has only finitely many nonzero rational homology groups, or dually, 
$Y$ has only finitely many nonzero rational homotopy groups \cite{MS}.
\item[(ii)] $Y$ is a bouquet of suspensions of connected finite complexes \cite{Gh, Le}. 
\item[(iii)] $\Ph(X,Y)$ has a natural abelian group structure and it has no elements of order $p$ 
for some prime $p$ \cite{Le}.  
\end{enumerate}
}

Next, we provide another collection of spaces for which the Gray index conjecture is valid. 
\\

\noindent
{\bf Theorem 1.5.}  {\it If $Y$ is a simply connected finite complex, then for 
any {\rm CW} complex $X$ of finite type we have $\Ph^{\infty}(X,Y)=*.$}
\\

The remainder of the paper is devoted to proving the foregoing results.  
Our proofs are based on the tower-theoretic approach to phantom maps. 

For an inverse tower $G=\{G_1 \leftarrow G_2  \leftarrow G_3 \leftarrow \cdots\}$ of groups, 
not necessarily abelian,  Bousfield and Kan \cite{BK} defined $\varprojlim^1G$ 
and proved that a short exact sequence of inverse towers 
\[
1\to K\to G\to H\to 1
\]
induces a six term $\varprojlim-\varprojlim^1$ exact sequence of pointed sets
\[
1\to \varprojlim K\to \varprojlim G \to  \varprojlim H \to 
 \varprojlim\!{}^1K\to \varprojlim\!{}^1G\to \varprojlim\!{}^1 H\to *. 
\]

They also established a short exact sequence of pointed sets
\[
*\to \lim_{\leftarrow n}\!{}^1[\Sigma X_n,Y]\to [X,Y]\to \lim_{\leftarrow n}[X_n,Y]\to *
\]
for any CW-complex $X$. By using this sequence and the definition of phantom maps, we obtain a bijection of pointed 
sets
\[
\Ph(X,Y)\cong \lim_{\leftarrow n}\!{}^1[\Sigma X_n,Y]. 
\]

\section{Proof of Theorem 1.2 and Corollary 1.3}

We start by constructing an algebraic example. 

Let 
\[
A=\Z/2^{\infty}\times \Z/3^{\infty}\times\cdots\times \Z/p^{\infty}\times\cdots,
\]
where $p$ is a prime and $\Z/p^{\infty}=\Q/\Z\p$. $A_n$ is a subgroup of $A$ defined by 
\[
A_n=\{(x_1,x_2,\cdots,x_n,x_{n+1},\cdots)\in A\ |\ x_1=x_2=\cdots=x_n=0\}.
\]
We embed $\Z$ in $A$ via the map 
\[
k\mapsto k(\frac{1}{2},\frac{1}{3},\cdots,\frac{1}{p},\cdots)
\]
and its image is also denoted by $\Z$. Now we consider the following tower of abelian groups
\[
A_0=A/\Z\leftarrow A_1\leftarrow A_2\leftarrow \cdots  \leftarrow A_n\leftarrow \cdots,
\]
where the first map is the composition of the inclusion and the projection and the others are all inclusion maps. 

For a tower $G=\{G_1 \leftarrow G_2  \leftarrow G_3 \leftarrow \cdots\}$ of sets, we define 
$G_n^{(k)}=\mathrm{Im}(G_k\to G_n)$ if $k\geq n$, and $G_n^{(k)}=G_k$ if $k<n$, 
and $G_n^{(\infty)}=\cap_{k\geq n}G_n^{(k)}$. 
\\

\noindent
{\bf Lemma 2.1.} {\it 
\[
A_0^{(\infty)}=(\Z/2\times \Z/3\times\cdots\times \Z/p\times\cdots)/\Z\cong \R\oplus\Q/\Z, 
\]
where $\R$ is viewed as a rational vector space whose cardinality equals that of the real numbers.}
\\

\noindent
{\bf Proof.} Let 
\[
[(\frac{n_2}{2},\frac{n_3}{3},\cdots,\frac{n_p}{p},\cdots)]\in (\Z/2\times \Z/3\times\cdots\times \Z/p\times\cdots)/\Z,
\]
then by the Chinese Remainder Theorem there is an integer $k$ such that 
$n_2\equiv k\pmod{2},n_3\equiv k\pmod{3},\cdots$ 
for the first $n$ primes. Thus 
\begin{eqnarray*}
&&[(\frac{n_2}{2},\frac{n_3}{3},\cdots,\frac{n_p}{p},\cdots)]\\
&=&[(\frac{n_2}{2},\frac{n_3}{3},\cdots,\frac{n_p}{p},\cdots)-k(\frac{1}{2},\frac{1}{3},\cdots,\frac{1}{p},\cdots)]\\
&=&[(0,\cdots,0,\frac{n_p-k}{p},\cdots)]\in A_0^{(n)},
\end{eqnarray*}
where $p$ is the $(n+1)$-st prime. As $n$ is arbitrary, 
we have $A_0^{(\infty)}\supset (\Z/2\times \Z/3\times\cdots\times \Z/p\times\cdots)/\Z$. 

To the contrary, let $[(q_2,q_3,\cdots,q_p,\cdots)]\in A_0^{(n)}$; 
then there is an integer $k$ and a set of rationals $\{r_p\}$ such that
\[
(q_2,q_3,\cdots,q_p,\cdots)=k(\frac{1}{2},\frac{1}{3},\cdots,\frac{1}{p},\cdots)+ (0,0,\cdots,0,r_p,\cdots),
\]
where $p$ is the $(n+1)$-st prime.
Thus $A_0^{(\infty)}\subset (\Z/2\times \Z/3\times\cdots\times \Z/p\times\cdots)/\Z$ 
and we proved the first assertion. 

As for the group structure, we make use of the computation of $\varprojlim^1$.  

It is easy to see that $K_n^0=\Ker(A_n\to A_0)=k(n) K\subset A_n\subset A$ for $n\geq1$, 
where $k(n)$ is a product of the first $n$ primes. 
Then there is a short exact sequence of towers:
\[
0\to \{K_n^0\}\to \{A_n\}\to \{A_0^{(n)}\}\to 0
\]
This induces the following six term $\varprojlim-\varprojlim^1$ exact sequence
\[
0\to \varprojlim_n K_n^0\to \varprojlim_n A_n \to  \varprojlim_n A_0^{(n)} \to 
 \varprojlim\!{}_n^1 K_n^0\to \varprojlim\!{}_n^1 A_n\to \varprojlim\!{}_n^1 A_0^{(n)}. 
\] 

As $\varprojlim_n A_n=\varprojlim\!{}_n^1 A_n=0$, we have an isomorphism
\[
A_0^{(\infty)}=\varprojlim_n A_0^{(n)}\cong \varprojlim\!{}_n^1 K_n^0\cong \varprojlim\!{}_n^1 k(n)\Z\cong 
\R\oplus\Q/\Z.
\]
As for the last equation see, for example, p.1228 of \cite{Mc}. \quad $\Box$
\\

Roitberg \cite{R1} constructs spaces $M^{2n+1}(I,J)$ to show 
that some groups $\Ph(X,Y)$ can possess torsion. Here 
we use these spaces to construct the space $Y$ stated in Theorem 1.2. 

Let $I$ and $J$ denote non-empty complementary sets of primes. The space 
 $M^{2n+1}(I,J)$, $n\geq2$, are defined by means of homotopy-pullback diagrams:
\[
\begin{CD}
M^{2n+1}(I,J)  @>>>  \Omega^2S_{J}^{2n+1}\\
@VVV                @VVV         \\
K(\Z_I,2n-1)   @>>>    K(\Q,2n-1)
\end{CD}
\] 
where $K(\Z_I,2n-1)$ and $S_J^{2n+1}$ denote the respective localizations of \linebreak 
$K(\Z,2n-1)$ and $S^{2n+1}$ and the bottom and right-hand maps are rationalizations. 

Let $p$ be a prime, $I_p=\{p\}$ and $J_p$ be the set of all primes except $p$. 
Set $M^{2p+3}=M^{2p+3}(I_p,J_p)$, that is,  $M^{2p+3}$ are defined by means of 
homotopy-pullback diagrams:
\[
\begin{CD}
M^{2p+3}=M^{2p+3}(I_p,J_p) @>>>  \Omega^2S_{J_p}^{2p+5}\\
@VVV                @VVV         \\
K(\Z_{I_p},2p+3)   @>>>    K(\Q,2p+3)
\end{CD}
\]
Note that $M^{2p+3}$ is a $2p+2$-connected double-loop space of finite type. 

By using a fiber sequence 
\[
K(\Q,2p+2)\overset{\delta}{\to} M^{2p+3} \to K(\Z_{I_p},2p+3)\times \Omega^2S_{J_p}^{2p+5},
\]
it is easy to show that 
\[
\Ph(\Sigma^2\mathbb{C}P^{\infty},M^{2p+3})=[\Sigma^2\mathbb{C}P^{\infty},M^{2p+3}]\cong \mathbb{R}\oplus \Z/p^{\infty}
\]
and that the torsion elements in $\Ph(\Sigma^2\mathbb{C}P^{\infty},M^{2p+3})$ 
come from \\ $[\Sigma^2\mathbb{C}P^{\infty},K(\Q,2p+2)]$. We use  
$\mathbb{C}P^{\infty}_m$ to denote the stunted projective space $\mathbb{C}P^{\infty}/\mathbb{C}P^{m-1}$. 
\\

\noindent
{\bf Lemma 2.2.} {\it If $p\geq m$, the canonical projection $\mathbb{C}P^{\infty}\to\mathbb{C}P^{\infty}_m$ induces an 
isomorphism 
\[
[\Sigma^2\mathbb{C}P^{\infty}_m, M^{2p+3}]=\Ph(\Sigma^2\mathbb{C}P^{\infty}_m, M^{2p+3})\to 
\Ph(\Sigma^2\mathbb{C}P^{\infty}, M^{2p+3}).
\]
If $p<m$, then $\Ph(\Sigma^2\mathbb{C}P^{\infty}_m, M^{2p+3})=0$.}
\\

\noindent
{\bf Proof.} As there is an epimorphism (see section 5 of \cite{Mc})
\begin{multline*}
\prod_kH^k(\Sigma^2\mathbb{C}P^{\infty}_m;\pi_{k+1}(M^{2p+3})\otimes\R)
\cong H^{2p+2}(\Sigma^2\mathbb{C}P^{\infty}_m;\R)\\
\to \Ph(\Sigma^2\mathbb{C}P^{\infty}_m, M^{2p+3}),
\end{multline*}
the assertion clearly follows for the case $p< m$. 

Next we consider the case $p\geq m$. We have 
\[
[\Sigma^2\mathbb{C}P^{\infty},M^{2p+3}]=\\
\Ph(\Sigma^2\mathbb{C}P^{\infty},M^{2p+3}). 
\]
Then the cofiber sequence $\mathbb{C}P^{m-1}\to \mathbb{C}P^{\infty}\to \mathbb{C}P^{\infty}_m\to\Sigma\mathbb{C}P^{m-1}$ 
induces an exact sequence 
\begin{multline*}
[\Sigma^3\mathbb{C}P^{m-1}, M^{2p+3}]\to [\Sigma^2\mathbb{C}P^{\infty}_m, M^{2p+3}]\\
\to\Ph(\Sigma^2\mathbb{C}P^{\infty}, M^{2p+3})\to[\Sigma^2\mathbb{C}P^{m-1}, M^{2p+3}].
\end{multline*}
As there is no essential phantom map in $[\Sigma^2\mathbb{C}P^{m-1}, M^{2p+3}]$ and \\
$[\Sigma^3\mathbb{C}P^{m-1}, M^{2p+3}]=0$ as $p\geq m$, this exact sequence induces an isomorphism:
\begin{equation}
[\Sigma^2\mathbb{C}P^{\infty}_m, M^{2p+3}]\cong\Ph(\Sigma^2\mathbb{C}P^{\infty}, M^{2p+3}).
\end{equation}

On the other hand, by Proposition 3 of \cite{LS}, we have 
\begin{multline}
\Ph(\Sigma^2\mathbb{C}P^{\infty}, M^{2p+3})=\Ph^{2m+2}(\Sigma^2\mathbb{C}P^{\infty}, M^{2p+3})\\
=\mathrm{Im}(\Ph(\Sigma^2\mathbb{C}P^{\infty}_m, M^{2p+3})\to \Ph(\Sigma^2\mathbb{C}P^{\infty}, M^{2p+3})).
\end{multline}

(2.1) and (2.2) imply the truth of the assertion for the case $p\geq m$. \quad $\Box$
\\

Let $\varphi_p$ be the composition 
\[
\varphi_p:\Sigma^2\mathbb{C}P^{\infty}\overset{\Sigma^2x^p/p}{\longrightarrow} \Sigma^2K(\Q,2p)
\overset{\iota}{\to} 
K(\Q,2p+2)\overset{\delta}{\to} M^{2p+3}
\]
where $x\in[\mathbb{C}P^{\infty}, K(\Z,2)]\cong H^2(\mathbb{C}P^{\infty};\Z)\subset H^2(\mathbb{C}P^{\infty};\Q)$ 
is the canonical generator and $\iota:\Sigma^2K(\Q,2p)\to K(\Q,2p+2)$ is the adjoint of the identity 
$K(\Q,2p)\to \Omega^2K(\Q,2p+2)$. $\varphi_p$ is a torsion element of order $p$. Now we put
\[
\varphi=(\prod_{p:\text{prime}} \varphi_p)\circ\Delta:\Sigma^2\mathbb{C}P^{\infty}\to 
\prod_{p:\text{prime}} \Sigma^2\mathbb{C}P^{\infty}\to \prod_{p:\text{prime}} M^{2p+3}
\]
and let $Y$ be the homotopy fiber of the map $\varphi$:
\[
Y\to \Sigma^2\mathbb{C}P^{\infty}\overset{\varphi}{\to} \prod_{p:\text{prime}} M^{2p+3}.
\]
Here $\prod_{p:\text{prime}} \Sigma^2\mathbb{C}P^{\infty}$ 
and $\prod_{p:\text{prime}} M^{2p+3}$ have the product topology. 

First we show that $Y$ is 3-connected and of finite type. By definition of $Y$ 
there is an exact sequence
\[
\pi_{i+1}( \prod_{p:\text{prime}} M^{2p+3})\to \pi_i(Y)\to \pi_i(\Sigma^2\mathbb{C}P^{\infty}).
\]
As $M^{2p+3}$ is $2p+2$-connected, we have an isomorphism 
\[
\pi_{i+1}( \prod_{p:\text{prime}} M^{2p+3})\cong \oplus_{2p+2\leq i}\pi_{i+1}(M^{2p+3}),
\]
therefore, it is finitely generated and 0 for $i<6$. 
Needless to say, $\pi_i(\Sigma^2\mathbb{C}P^{\infty})$ is finitely generated and 0 for $i<4$. 
By the exact sequence above $\pi_i(Y)$ is also finitely generated and 0 for $i<4$.

 Again by definition, there is an exact sequence 
\[
[\mathbb{C}P^{\infty}_m, \Omega^2\Sigma^2\mathbb{C}P^{\infty}]\overset{(\Omega^2 \varphi)_*}{\to} 
[\mathbb{C}P^{\infty}_m, \prod_{p:\text{prime}} \Omega^2 M^{2p+3}]\to [\mathbb{C}P^{\infty}_m,\Omega Y].
\]
Thus there are inclusions 
\[
\Ph(\mathbb{C}P^{\infty}_m, \prod_{p:\text{prime}} \Omega^2 M^{2p+3})/\mathrm{Im}(\Omega^2\varphi)_*\subset 
\Ph(\mathbb{C}P^{\infty}_m,\Omega Y)
\]
and we analyze the tower 
\begin{multline*}
\{\Ph(\mathbb{C}P^{\infty}_m, \prod_{p:\text{prime}} \Omega^2 M^{2p+3})/\mathrm{Im}(\Omega^2 \varphi)_*\}_m
\\
\cong \{\Ph(\Sigma^2\mathbb{C}P^{\infty}_m, \prod_{p:\text{prime}} M^{2p+3})/\mathrm{Im}\varphi_*\}_m.
\end{multline*}
\\

\noindent
{\bf Proposition 2.3.} {\it The map
\[
(\varphi_p)_*:[\Sigma^2\mathbb{C}P^{\infty}_m,\Sigma^2\mathbb{C}P^{\infty}]\to 
[\Sigma^2\mathbb{C}P^{\infty}_m,M^{2p+3}].
\]
is trivial for $m>1$. For $m=1$, its image is isomorphic to $\Z/p$ and is generated by $\varphi_p$.}
\\

\noindent
{\bf Proof.} As $\varphi_p$ factors through $K(\Q,2p+2)$, it is sufficient to establish a set of 
generators of $[\Sigma^2\mathbb{C}P^{\infty}_m,\Sigma^2\mathbb{C}P^{\infty}]$ 
up to homology to compute $\mathrm{Im}\ \varphi_*$. 
McGibbon \cite{Mc82} proves that any self-map of $\Sigma^k\mathbb{C}P^{\infty}$ 
is homologous to a linear combination of $k$-fold suspensions of elements in 
$[\mathbb{C}P^{\infty},\mathbb{C}P^{\infty}]\cong\Z$ for $k\geq1$. 
Another set of generators is given by Morisugi \cite{Mo}. 

Let $f_1^{\top}:\mathbb{C}P^{\infty}\to \Omega\Sigma\mathbb{C}P^{\infty}$ be the adjoint of the identity 
$f_1:\Sigma\mathbb{C}P^{\infty}\to \Sigma\mathbb{C}P^{\infty}$. 
Following \cite{Mo}, we define inductively  
\[
f_{n+1}^{\top}:\mathbb{C}P^{\infty}\overset{\tilde{\Delta}}{\to}\mathbb{C}P^{\infty}\wedge\mathbb{C}P^{\infty}
\overset{f_1^{\top}\wedge f_n^{\top}}{\to}\Omega\Sigma\mathbb{C}P^{\infty}\wedge\Omega\Sigma\mathbb{C}P^{\infty}
\overset{\sharp}{\to}\Omega\Sigma\mathbb{C}P^{\infty},
\]
where $\tilde{\Delta}$ is the reduced diagonal map and $\sharp$ denotes an extension of the adjoint of the Hopf 
construction of $\mathbb{C}P^{\infty}$. Morisugi proves that $f_n^{\top}$ factors as 
\[
\mathbb{C}P^{\infty}\to \mathbb{C}P^{\infty}_n=\mathbb{C}P^{\infty}/\mathbb{C}P^{n-1}\overset{g_n^{\top}}{\to} 
\Omega\Sigma\mathbb{C}P^{\infty},
\]
where the first map is the canonical projection. Let $f_n:\Sigma\mathbb{C}P^{\infty}\to\Sigma\mathbb{C}P^{\infty}$ 
and $g_n:\Sigma\mathbb{C}P^{\infty}_n\to \Sigma\mathbb{C}P^{\infty}$ be the adjoint of $f_n^{\top}$ and 
$g_n^{\top}$, respectively. $\{\Sigma^{k-1}f_n\}_{n\geq1}$ is a set of generators of self-maps of 
$\Sigma^k\mathbb{C}P^{\infty}$ up to homology for $k\geq1$.

Let $\beta_n\in H_{2n}(\mathbb{C}P^{\infty};\Z)\cong\Z$ be the dual of $x^n$, 
where  $H^*(\mathbb{C}P^{\infty};\Z)=\Z[x]$. We put  
\[
(f_n)_*(\sigma\beta_k)=\delta_n(k)\sigma\beta_k,
\]
where $\sigma:\tilde{H}_*(\mathbb{C}P^{\infty};\Z)\to \tilde{H}_{*+1}(\Sigma\mathbb{C}P^{\infty};\Z)$ 
is the suspension isomorphism. Then $\delta_n(k)$ is given by the following formula (see, Theorem 3.3 of \cite{Mc82}), 
\[
\delta_n(k)=\sum_{i=1}^n(-1)^{n-i}\binom{n}{i}i^k.
\]
$\delta_n(k)$ is known to be 0 for $k<n$ and to be divisible by $n!$. 
\\

\noindent
{\bf Lemma 2.4.} {\it If $p$ is a prime, then $\delta_n(p)$ is divisible by $p$ for $n>1$.}
\\

\noindent
{\bf Proof.} As $i^p\equiv i\pmod{p}$, we have 
\begin{eqnarray*} 
\delta_n(p)&=&\sum_{i=1}^n(-1)^{n-i}\binom{n}{i}i^p\\
           &\equiv&\sum_{i=1}^n(-1)^{n-i}\binom{n}{i}i\pmod{p}\\
           &\equiv&\delta_n(1)\pmod{p}\\
           &\equiv& 0\pmod{p}.
\end{eqnarray*}
\\

It is easy to see that 
$\{\Sigma^2p_k^m\circ \Sigma g_k:\Sigma^2\mathbb{C}P^{\infty}_m\to \Sigma^2\mathbb{C}P^{\infty}|k\geq m\}$ is 
a set of generators of $[\Sigma^2\mathbb{C}P^{\infty}_m, \Sigma^2\mathbb{C}P^{\infty}]$ up to homology, where 
$p^m_k:\mathbb{C}P^{\infty}_m\to \mathbb{C}P^{\infty}_k$ is the canonical projection for $k\geq m$. 
Consider the following commutative diagram 
\[
\begin{CD}
[\Sigma^2 \mathbb{C}P^{\infty}_m, \Sigma^2 \mathbb{C}P^{\infty}]@>(\varphi_p)_*>> \Ph(\Sigma^2 \mathbb{C}P^{\infty}_m,M^{2p+1})\\
@VVj_1V                        @VVj_2V \\
[\Sigma^2 \mathbb{C}P^{\infty}, \Sigma^2 \mathbb{C}P^{\infty}]@>(\varphi_p)_*>> \Ph(\Sigma^2 \mathbb{C}P^{\infty},M^{2p+1})\\
\end{CD}
\]
where the vertical maps are induced by the projection $\Sigma^2 \mathbb{C}P^{\infty}\to \Sigma^2 \mathbb{C}P^{\infty}_m$. 
As $j_1$ is injective up to homology, $j_2$ is injective and $\Sigma f_1$ is the identity, 
to prove Proposition 2.3, it is sufficient to prove the following.  
\\

\noindent
{\bf Lemma 2.5.} {\it $(\varphi_p)_*(\Sigma f_n)=0$ for $n>1$.}
\\

\noindent
{\bf Proof.} As
\[
(\Sigma^2\frac{x^p}{p}\circ\Sigma f_n)_*(\sigma^2(\beta_p))=\delta_n(p)(\Sigma^2\frac{x^p}{p})_*(\sigma^2(\beta_p)),
\]
we have $(\varphi_p)_*(\Sigma f_n)=\delta_n(p)\varphi_p=0$. 
Here we use the fact that $\varphi_p$ is of order $p$ and 
$\delta_n(p)$ is divisible by $p$, due to Lemma 2.4.  \quad $\Box$ 
\\

\noindent
{\bf Proof of Theorem 1.2.} We have inclusions 
\[
\Ph(\Sigma^2\mathbb{C}P^{\infty}_m, \prod_{p:\text{prime}} M^{2p+3})/\mathrm{Im}(\varphi)_*\subset 
\Ph(\mathbb{C}P^{\infty}_m,\Omega Y) 
\]
and in the tower 
\[
\{\Ph(\Sigma^2\mathbb{C}P^{\infty}_m, \prod_{p:\text{prime}} M^{2p+3})/\mathrm{Im}\varphi_*\}_m,
\]
we have 
\begin{multline*}
\Ph(\Sigma^2\mathbb{C}P^{\infty}, \prod_{p:\text{prime}} M^{2p+3})/\mathrm{Im}\varphi_*\\
\cong ((\Z/2^{\infty}\oplus\R)\times (\Z/3^{\infty}\oplus\R)\times\cdots\times(\Z/p^{\infty}\oplus\R)\times\cdots)/K
\end{multline*}
where $K\cong\Z$ is a subgroup generated by $(\varphi_p)_p$, and for $m>1$ we have 
\begin{eqnarray*}
\Ph(\Sigma^2\mathbb{C}P^{\infty}_m, \prod_{p:\text{prime}} M^{2p+3})/\mathrm{Im}\varphi_*
&=&\Ph(\Sigma^2\mathbb{C}P^{\infty}_m, \prod_{p:\text{prime}} M^{2p+3})\\
&\cong& (\Z/\ell^{\infty}\oplus\R)\times\cdots\times(\Z/p^{\infty}\oplus\R)\times\cdots
\end{eqnarray*}
where $\ell$ is the smallest prime which is equal to $m$ or larger than $m$. Then Lemma 2.1 shows that 
\[
\Ph^{\infty}(\mathbb{C}P^{\infty},\Omega Y)\supset 
(\Ph(\Sigma^2\mathbb{C}P^{\infty}, \prod_{p:\text{prime}} M^{2p+3})/\mathrm{Im}(\varphi)_*)^{(\infty)}\neq0
\]
and the proof is complete. \quad $\Box$ 
\\

We recall the natural filtration on $\SNT(X)$ according to Ghienne \cite{Gh}. 

Let $\{G_n\}_n$ be an inverse tower of groups. 
A surjection of towers $\{G_n\}_n\to \{G_k^{(n)}\}_n$ induces a surjection of $\varprojlim^1$ sets: 
\[
p_k:\lim_{\leftarrow n}\!{}^1G_n\to \lim_{\leftarrow n}\!{}^1G_k^{(n)}
\]
Set $L=\displaystyle{\lim_{\leftarrow n}\!{}^1G_n}$ and define $L^k=\Ker\ p_k$. We then have 
a filtration:
\[
L=L^0\supset L^1\supset \cdots \supset L^k\supset\cdots\supset L^{\infty}=\cap_kL^k,
\]
which is called the algebraic Gray filtration on $L=\displaystyle{\lim_{\leftarrow n}\!{}^1G_n}$. 

For a connected space $X$, we denote by $\mathrm{Aut}(X)$ the group of 
homotopy classes of self-homotopy equivalences of $X$. Recall from \cite{W} 
that we have a bijection:
\[
\SNT(X)\cong \lim_{\leftarrow n}\!{}^1\mathrm{Aut}(X^{(n)})
\]
This description of $\SNT(X)$ as $\varprojlim^1$ set defines the algebraic Gray filtration on it:
\[
\SNT(X)\supset \SNT^1(X)\supset \cdots \supset \SNT^k(X)\supset\cdots
\supset \SNT^{\infty}(X)=\cap_k\SNT^k(X)
\]

\noindent
{\bf Proof of Corollary 1.3.} For a phantom map $f:X\to Z$, its homotopy fiber has the same 
$n$-type as $X\times\Omega Z$ for all $n$ since $f^{(n)}:X^{(n)}\to Y^{(n)}$ is null-homotopic. 
Then we can define a map
\[
\mathrm{F}:\Ph(X,Z)\to \SNT(X\times\Omega Z)
\]
which associates to a phantom map its homotopy fiber. 
Theorem 3.6 of \cite{Gh} says that this map respects 
filtration. In particular, $\mathrm{F}(\Ph^{\infty}(X,Z))\subset  \SNT^{\infty}(X\times\Omega Z)$. 

To prove Corollary 1.3, therefore, it is sufficient to prove that there is a phantom map 
$\phi\in \Ph^{\infty}(\mathbb{C}P^{\infty},\Omega Y)$ such that its homotopy fiber $F_{\phi}$ is not homotopy equivalent to 
$\mathbb{C}P^{\infty}\times\Omega^2 Y$. 

We follow the argument of the proof of Lemma 3.3 in \cite{HR}. 

Let $\phi$ be an element of $\Ph^{\infty}(\mathbb{C}P^{\infty},\Omega Y)$ with infinite order which factors through $\prod\Omega M^{2p+3}$. 

Let $f:\mathbb{C}P^{\infty}\times \Omega^2Y\to F_{\phi}$ be any map, $j:\mathbb{C}P^{\infty}\to 
\mathbb{C}P^{\infty}\times \Omega^2Y$ the canonical embedding and consider the fiber sequence 
\[
\to \Omega^2Y\overset{b}{\to}F_{\phi}\overset{i}{\to}\mathbb{C}P^{\infty}\overset{\phi}{\to}\Omega Y.
\]
We write $d=i\circ f\circ j\in[\mathbb{C}P^{\infty},\mathbb{C}P^{\infty}]\cong\Z$. 
We have $\phi\circ d=\phi\circ i\circ f\circ j=0$ as $\phi\circ i=0$.  

By Lemma 2.1 we can identify $\phi$ with an element in 
\[
(\Z/2\times \Z/3\times\cdots\times \Z/p\times\cdots)/\Z
\]
and write 
\[
\phi=[(\frac{n_2}{2},\frac{n_3}{3},\cdots,\frac{n_p}{p},\cdots))].
\]
Then we can write 
\begin{align*}
\phi\circ d&=[(d^2\frac{n_2}{2},d^3\frac{n_3}{3},\cdots,d^p\frac{n_p}{p},\cdots)],\\
\intertext{and as $d^p\equiv d\pmod{p}$, we have}
           &=[(d\frac{n_2}{2},d\frac{n_3}{3},\cdots,d\frac{n_p}{p},\cdots)]=d\phi.
\end{align*}
As $\phi$ is an element of infinite order, $\phi\circ d=d\phi=0$ implies that $d=0$. 
 
Then there is a map $g:\mathbb{C}P^{\infty}\to \Omega^2Y$ such that $f\circ j\simeq b\circ g$. 

Consider the following commutative diagram:
\[
\begin{CD}
[\mathbb{C}P^{\infty},\Omega^2Y]@>\xi_*>>[\mathbb{C}P^{\infty},\Omega^2\Sigma^2\mathbb{C}P^{\infty}]
@>(\Omega^2\varphi)_*>>[\mathbb{C}P^{\infty},\prod\Omega^2M^{2p+3}]\\
@VV\ell^*V           @VV\ell^*V           @VV\ell^*V \\
\pi_2(\Omega^2Y)@>\xi_*>>\pi_2(\Omega^2\Sigma^2\mathbb{C}P^{\infty})
@>(\Omega^2\varphi)_*>>\pi_2(\prod\Omega^2M^{2p+3})
\end{CD}
\]
Here the horizontal sequences are induced by the fiber sequence $ \Omega^2Y\to \Omega^2\Sigma^2\mathbb{C}P^{\infty}
\to \prod\Omega^2M^{2p+3}$ and vertical maps are induced by the canonical inclusion 
$\ell:S^2=\mathbb{C}P^1\to \mathbb{C}P^{\infty}$. In the proof of Theorem 1.2 we proved that the kernel of the map 
$(\Omega^2\varphi)_*:[\mathbb{C}P^{\infty},\Omega^2\Sigma^2\mathbb{C}P^{\infty}]\to 
[\mathbb{C}P^{\infty},\prod\Omega^2M^{2p+3}]$ is generated by the maps which are the adjoints of $\{\Sigma f_k\ |\ k\geq2\}$ 
up to homology. Thus, as $\xi_*(g)$ factors through $\mathbb{C}P^{\infty}_2$ up to homology, we have 
\[
0=\ell^*\circ\xi_*(g)=\xi_*\circ\ell^*(g),
\]
that is, 
\begin{align*}
\ell^*(g)&\in \mathrm{Ker}(\xi_*:\pi_2(\Omega^2Y)\to \pi_2(\Omega^2\Sigma^2\mathbb{C}P^{\infty}))\\
          &\in \mathrm{Im}(\pi_2(\prod\Omega^3M^{2p+3})\to \pi_2(\Omega^2Y)).
\end{align*}
On the other hand, $\pi_2(\prod\Omega^3M^{2p+3})\cong \prod\pi_5(M^{2p+3})=0$ as $M^{2p+3}$ is $2p+2$-connected. 
It follows that 
\[
g_*:\pi_2(\mathbb{C}P^{\infty})\to \pi_2(\Omega^2Y)
\]
is the 0 map. As $f_*\circ j_*=b_*\circ g_*=0$ on 2-dimensional homotopy groups and 
$j_*:\pi_2(\mathbb{C}P^{\infty})\to \pi_2(\mathbb{C}P^{\infty}\times\Omega^2 Y)$ is an embedding, we see that 
\[
f_*:\pi_2(\mathbb{C}P^{\infty}\times\Omega^2 Y)\to \pi_2(F_{\phi})
\]
maps the factor $\pi_2(\mathbb{C}P^{\infty})$ trivially, hence, that $f$ cannot be a homotopy equivalence. 
\quad $\Box$

\section{Proof of Theorem 1.5.}

By using a bijection of pointed sets given by L\^{e} Minh H\`{a} \cite{Le} 
\[
\Ph^{\infty}(X,Y)\cong \lim_{\leftarrow n}\!{}^1[\Sigma X_n,Y]^{(\infty)}, 
\]
it is sufficient to prove that the tower $\{[\Sigma X_n,Y]^{(\infty)}\}_n$ satisfies 
the Mittag-Leffler condition to prove Theorem 1.5. 
From now on in this proof, let $G_n=[\Sigma X_n,Y]\cong[X_n,\Omega Y]$. 

To prove that the tower $\{G_n^{(\infty)}\}_n$ satisfies the Mittag-Leffler condition, 
it is sufficient to prove that 
for each $n$, the image of the map 
\[
\lim_{\leftarrow k}G_k=\lim_{\leftarrow k}G_k^{(\infty)}\to G_n^{(\infty)}
\]
has a finite cokernel. For $\mathrm{Im}(G_m^{(\infty)}\to G_n^{(\infty)})$ has only finitely many 
possibilities in $G_n^{(\infty)}$ as  
\[
\mathrm{Im}(\lim_{\leftarrow k}G_k^{(\infty)}\to G_n^{(\infty)})\subset 
\mathrm{Im}(G_m^{(\infty)}\to G_n^{(\infty)})\subset G_n^{(\infty)}.
\]

Since there is a surjection $[X,\Omega Y]\to \lim_{\leftarrow k}G_k$, therefore, 
\[
\mathrm{Im}([X,\Omega Y]\to G_n^{(\infty)})=\mathrm{Im}(\lim_{\leftarrow k}G_k\to G_n^{(\infty)}),
\] 
it is sufficient to show that the image $[X,\Omega Y]\to G_n^{(\infty)}$ has a finite cokernel. 

When $\Omega Y$ has the rational homotopy type of 
\[
P_1=\prod_{\alpha\in A_1} S^{2n_{\alpha}+1}\times \prod_{\alpha\in A_2} \Omega S^{2n_{\alpha}+1}, 
\]
where $P_1$ is topologized as the direct limit of finite products, 
it is easy to construct a rational homotopy equivalence
\[
f:P_1\to \Omega Y.
\] 
Here a map $f:X\to Y$ is called a rational homotopy equivalence if it induces 
an isomorphism in rational homology groups. 

In \cite{I} we construct a rational homotopy equivalence in the opposite direction. 
We modify this result to prove Theorem 1.5 as follows. 

\bigskip

\noindent
{\bf Lemma 3.1.}  {\it For a natural number $n$ and a rational homotopy equivalence $g:\Omega Y\to P_1$ 
there is a rational homotopy equivalence $f:P_1\to \Omega Y$ such that 
\[
(f\circ g)^{(n)}:(\Omega Y)^{(n)}=\Omega Y^{(n+1)}\to (\Omega Y)^{(n)}=\Omega Y^{(n+1)}
\]
is a power map.}
\\

Assume for the moment that this lemma is true and we continue the proof of Theorem 1.5. 

We set 
\[
P_2=\prod_{\alpha\in A_1} \Omega S^{2n_{\alpha}+2}\times \prod_{\alpha\in A_2} \Omega S^{2n_{\alpha}+1}=
\prod_{\alpha\in A}\Omega S^{n_{\alpha}}.
\]
Let $\varphi=\prod_{\alpha\in A_1} i_{\alpha}\times id_{\prod_{\alpha\in A_2} \Omega S^{2n_{\alpha}+1}}:P_1\to P_2$, 
where $i_{\alpha}: S^{2n_{\alpha}+1}\to \Omega S^{2n_{\alpha}+2}$ is the adjoint of the identity of $S^{2n_{\alpha}+2}$. 
Any map $h:S^m\to \Omega Y$ is equal to the composition $\Omega h^{\top}\circ i:S^m\to \Omega S^{m+1}\to \Omega Y$, 
where $i:S^m\to \Omega S^{m+1}$ and $h^{\top}:S^{m+1}\to Y$ are the adjoint of $id:S^{m+1}\to S^{m+1}$ and 
$h:S^m\to \Omega Y$, respectively. It follows that $f:P_1\to \Omega Y$ is factored as 
\[
f=f^{\prime}\circ \varphi:P_1\to P_2\to \Omega Y. 
\]

Now consider the following commutative diagram:
\begin{equation}
\begin{CD}
[X,\Omega Y]@>g_*>> [X, P_1]@>\varphi_*>> [X, P_2]@>f^{\prime}_*>>[X,\Omega Y]\\
  @VVV       @VVV    @VVV    @VVV\\
[X_n,\Omega Y]^{(\infty)}@>g_*>>[X_n, P_1]^{(\infty)}@>\varphi_*>> 
[X_n, P_2]^{(\infty)}@>f^{\prime}_*>>[X_n,\Omega Y]^{(\infty)}\\
\end{CD}
\end{equation}
where the vertical maps are induced by the inclusion map $X_n\to X$. 

First we assert that the second vertical map from the right is surjective.  
To prove this it is sufficient to prove that for each $m$, 
we have $\mathrm{Im}([X,\Omega S^m]\to [X_n,\Omega S^m])=[X_n,\Omega S^m]^{(\infty)}$. According to the remark after 
Proposition 2.6 of \cite{Le}, $\{[X_n,\Omega S^m]^{(\infty)}\}_n$ is a tower of epimorphisms. 
Then all maps in the sequence 
\[
[X,\Omega S^m]\to \lim_{\leftarrow n}[X_n,\Omega S^m]
=\lim_{\leftarrow n}[X_n,\Omega S^m]^{(\infty)}\to [X_n,\Omega S^m]^{(\infty)}
\]
are epimorphic and we complete the proof of the assertion. 

For each $n$, by Lemma 3.1, there are rational homotopy equivalences $f:P_1\to \Omega Y$ and 
$g:\Omega Y\to P_1$ such that 
\[
(f\circ g)^{(n)}:(\Omega Y)^{(n)}=\Omega Y^{(n+1)}\to (\Omega Y)^{(n)}=\Omega Y^{(n+1)}
\]
is a power map, say $\lambda$. By $\lambda$ we denote both a natural number and the power map on 
$\Omega Y$ of power $\lambda$. 
Then for each $x\in[X_n,\Omega Y]^{(\infty)}$ there is 
$u\in [X_n,\prod_{\alpha}\Omega S^{n_{\alpha}}]^{(\infty)}$ 
such that $f^{\prime}_*(u)=x^{\lambda}$. In fact, put $u=\varphi_*\circ g_*(x)$, 
then we have $f^{\prime}_*(u)=f^{\prime}_*(\varphi_*\circ g_*(x))=(f\circ g)_*(x)=
(f\circ g)^{(n)}_*(x)=\lambda_*(x)=x^{\lambda}$. Then, 
by the commutative diagram (3.1) and the fact that $[X,P_2]\to [X_n,P_2]^{(\infty)}$ is surjective, 
we conclude that for each $x\in[X_n,\Omega Y]^{(\infty)}$ there exist $v\in [X,\Omega Y]$ and 
a natural number $\lambda$ such that $v$ is mapped to $x^{\lambda}$. Thus by Lemma 7.1.2 of \cite{Mc} we conclude that  
the image $[X,\Omega Y]\to [X_n,\Omega Y]^{(\infty)}$ has finite index in $[X_n,\Omega Y]^{(\infty)}$ and 
complete the proof of Theorem 1.5. \quad $\Box$ 
\\

\noindent
{\bf Proof of Lemma 3.1.} For a nilpotent space $X$ by $r=r(X):X\to X_{(0)}$ we denote the rationalization of $X$. 
Consider the following commutative diagram:
\[
\begin{CD}
\Omega Y @>r>>      (\Omega Y)_{(0)}\\
@VVg V                  @VVg_{(0)}V  \\
P_1     @>r>>        (P_1)_{(0)}
\end{CD}
\]
As $g_{(0)}$ is a homotopy equivalence, there is a homotopy inverse $f^{\prime\prime}:(P_1)_{(0)}\to  (\Omega Y)_{(0)}$. 
Thus $r(\Omega Y)\simeq f^{\prime\prime}\circ r(P_1)\circ g=f^{\prime}\circ g$, 
where $f^{\prime}=f^{\prime\prime}\circ r(P_1):P_1\to (\Omega Y)_{(0)}$. 
  
Put 
\[
P^{\prime}=\prod_{\alpha\in A_1,2n_{\alpha}+1\leq n+1}S^{2n_{\alpha}+1}\times 
\prod_{\alpha\in A_2,2n_{\alpha}\leq n+1}(\Omega S^{2n_{\alpha}+1})_{n+1},
\]
where $(\Omega S^{2n_{\alpha}+1})_{n+1}$ denotes the $n+1$-skeleton of $\Omega S^{2n_{\alpha}+1}$.
As the rationalization of $\Omega Y$ can be constructed as an infinite telescope using power maps and we may assume that 
$g$ is a cellular map, there is a sufficiently large power map $\lambda:\Omega Y\to \Omega Y$ such that
\[
\lambda|_{(\Omega Y)_{n+1}}\simeq h\circ g|_{(\Omega Y)_{n+1}}.
\]
Here $h:P^{\prime}\to \Omega Y$ and $f^{\prime}|_{P^{\prime}}\simeq j\circ h$ where 
$j:\Omega Y\to  (\Omega Y)_{(0)}$ is an inclusion to the telescope. Then we have 
\[
\lambda^{(n)}\simeq h^{(n)}\circ g^{(n)}.
\]
As 
\[
[\Omega S^{2m+1}, \Omega Y]\cong [\Sigma\Omega S^{2m+1}, Y]\cong \prod_k[S^{2mk+1}, Y],
\]
every map $\ell^{\prime}:(\Omega S^{2m+1})_{n+1}\to \Omega Y$ has an extension 
$\ell:\Omega S^{2m+1}\to \Omega Y$.  Thus it is easy to construct a rational 
equivalence $f:P_1\to \Omega Y$ such that $f^{(n)}\simeq h^{(n)}$. 
Then $f$ satisfies the required condition in Lemma 3.1. 
\quad $\Box$

\bigskip

\noindent
{\bf Acknowledgments}

\bigskip

The author thanks D. Kishimoto for his valuable comments on an earlier draft. 



\end{document}